\documentclass[11pt]{amsart}
\usepackage{amssymb,latexsym,comment,url,amsmath}

\usepackage[T1]{fontenc}
\usepackage{rotating,graphicx}
\usepackage{xcolor}
\usepackage{diagbox} 
\usepackage{soul}

\usepackage{longtable}

\newcommand{\End}{{\operatorname{End}}}
\newcommand{\Aut}{{\operatorname{Aut}}}

\newcommand{\Char}{\operatorname{char}}
\newcommand{\Jac}{\operatorname{Jac}}

\newcommand{\isom}{ \cong }

\newcommand{\OK}{{\mathcal{O}_K}}

\newcommand{\Q}{{\mathbb Q}}

\newcommand{\Disc}{\operatorname{Discriminant}}

\newcommand{\Z}{{\mathbb Z}}

\newcommand{\Magma}{{\sf MAGMA }}
\newcommand{\Sage}{{\sf Sage }}
\newcommand{\Mathematica}{{\sf Mathematica }}

\newenvironment{Proof}{\par\noindent{\sc Proof:}}%
                      {\hspace*{\fill}\nobreak$\Box$\par\medskip}
                       {\hspace*{\fill}\nobreak$\Box$\par\medskip}

\newtheorem{Proposition}{Proposition}[section]
\newtheorem{Theorem}[Proposition]{Theorem}

\newtheorem{Algorithm}[Proposition]{Algorithm}
\newtheorem{Corollary}[Proposition]{Corollary}

\theoremstyle{definition}

\newtheorem{Remark}[Proposition]{Remark}
 \newtheorem{Example}[Proposition]{Example}

\renewcommand{\baselinestretch}{1.1}

\begin{document}

\title[Genus two curves with everywhere good reduction over quadratic fields]%
{Genus two curves with everywhere good reduction over quadratic fields}

\author[A. D\k{a}browski]%
{Andrzej~D\k{a}browski}
\address{Institute of Mathematics, University of Szczecin, Wielkopolska 15, 70-451 Szczecin, Poland}
\email{dabrowskiandrzej7@gmail.com}
\email{andrzej.dabrowski@usz.edu.pl}

\author[M. Sadek]%
{Mohammad~Sadek}
\address{Faculty of Engineering and Natural Sciences, Sabanc{\i} University, Tuzla, \.{I}stanbul, 34956 Turkey}
\email{mohammad.sadek@sabanciuniv.edu}

\begin{abstract}
{We address the question of existence of absolutely simple abelian varieties of dimension 2 with everywhere good reduction over quadratic fields.  The emphasis will be given to the construction of pairs $(K,C)$, where $K$ is a quadratic number field and $C$ is a genus $2$ curve with everywhere good reduction over $K$. We provide the first infinite sequence of pairs $(K,C)$ where $K$ is a real (respectively, complex) quadratic field and $C$ has everywhere good reduction over $K$. Moreover, we show that the Jacobian of the curves $C$ in these families is an absolutely simple abelian variety.}
\end{abstract}

\maketitle

\let\thefootnote\relax

\section{Introduction}

It is due to Fontaine, \cite{Fontaine} (see also Abrashkin \cite{Abrashkin}), that there are no abelian varieties with everywhere 
good reduction over any of the number fields $\Q$, $\Q(\sqrt{-1})$, $\Q(\sqrt{-3})$, or $\Q(\sqrt{5})$. 
The first example of an elliptic curve with everywhere good reduction over a quadratic field
was given by Tate over the field $\mathbb Q(\sqrt{29})$, namely, the elliptic curve described by the Weierstrass equation $y^2 + xy +\left({5+\sqrt{29} \over 2}\right)^2 y = x^3$. Much work has been directed since then to giving examples of elliptic curves with everywhere good reduction over number fields with due attention to quadratic number fields, see 
\cite{Kagawa,Setzer,Stroeker,Zhao} and references therein. 
In \cite{CL}, an algorithm was given
to find all elliptic curves over a number field with good reduction outside a given
finite set of primes.

Shafarevich conjectured \cite{Sha} that for each number field $K$, finite set of places $S$, and an integer $g\ge 2$, there are only finitely many $K$-isomorphism classes of curves of genus $g$ over $K$ with good reduction outside $S$. The proof was sketched by him in the hyperelliptic case; for details see  \cite{Oort,Parshin}. According to Fontaine's result, if $K=\Q$ then $S$ must contain at least one prime so that such curves exist.  
 The papers \cite{DS,MS,Row} treated genus two curves over $\Q$ with bad reduction
at exactly one prime. 

Fixing integers $g,k \geq 2$, one may pose the following question:
does there exist a genus $g$ curve with everywhere good reduction over
some number field $K$ of degree $k$?  Note that if a curve $C$ has
everywhere good reduction over some number field, then the jacobian of
$C$ also has such a property, but the converse is not true in general.  
In \cite{BaazizBoxall}, examples of genus 2 curves, hence abelian varieties of dimension 2, 
 with everywhere good reduction over quartic fields were given.  However,  the Jacobian varieties 
of these curves are not absolutely simple. The first examples of genus two curves with everywhere good reduction
over real quadratic fields were constructed in 
\cite{DembeleKumar,Dembele}.  
This was approached by combining the computation of Hecke eigenvalues of certain Hilbert modular forms, using new rational models of Hilbert modular surfaces, together with assuming the Eichler-Shimura conjecture for Hilbert modular forms. Up to the knowledge of the authors, the only example of a genus 2 curve with everywhere good reduction over a complex quadratic field can be found in \cite[Theorem 5.2]{Berger}.

In this article we produce genus two curves with everywhere good reduction over real quadratic fields without the current restrictions on the discriminant or the class number of the field. We also display a $1$-parameter family of pairs of genus 2 curves and real (respectively, complex) quadratic fields such that the curve has everywhere good reduction over the quadratic field.  Finally we prove the existence of genus two curves with everywhere good reduction over a quadratic field such that its Jacobian has trivial endomorphisms, hence the Jacobian is an absolutely simple abelian variety. 

Our approach makes use of Mestre's algorithm to construct genus 2 curves from their Igusa invariants. Starting from rational invariants that satisfy certain local properties, we obtain genus two curves over the rational field. The only nontrivial automorphism of these curves is the hyperelliptic involution. Moreover, the minimal discriminants of these curves are of prescribed valuations at every rational bad prime of the curve. After a quadratic base change we prove that these curves attain everywhere good reduction. Mestre's algorithm can be adjusted to produce genus 2 curves with non-generic automorphism groups, \cite{CardonaQuer}. This enables us to produce infinitely many such curves with everywhere good reduction over corresponding quadratic fields.  
 \section*{Acknowledgement}
 The authors are very grateful to the anonymous referees for many corrections and valuable suggestions that improved the manuscript. All the calculations in this work were performed using \Magma \cite{Magma} and \Mathematica \cite{Mathematica}. M. Sadek is partially supported by BAGEP Award of the Science Academy, Turkey, and The Scientific and Technological Research Council of Turkey, T\"{U}B\.{I}TAK; research grant: ARDEB 1001/122F312.

\section{Preliminaries on genus two curves }
\label{sec:1}
 Let $K$ be a perfect field.
 Let $C$ be a (projective smooth) curve of genus two defined by a 
{\em Weierstrass equation} over $K$ of the form 
 \begin{eqnarray}\label{eq1} E:\;y^2+Q(x)y=P(x),\qquad P(x),Q(x)\in K[x],\; \deg Q(x)\le 3,\; 
 \deg P(x)\le 6.\end{eqnarray} 
 Given a Weierstrass equation of the form 
\begin{equation} \label{eqqqq}
y^2=f(x)=a_0x^6+a_1x^5+a_2x^2+a_3x^3+a_4x^2+a_5x+a_6,
\end{equation}
 Igusa associated to $f(x)$ the invariants $J_2$, $J_4$, $J_6$, $J_{10}\in \Z[1/2,a_0,a_2,\cdots,a_6].$ In fact, $J_{10}=2^{-12}\cdot\Disc(f)$ if $a_0\ne 0$, and $J_{10}=2^{-12}a_1^2\cdot\Disc(f)$ otherwise. For the full description of these invariants, see \cite{Igusa} when $a_0=0$, and \cite[\S2]{Liu3} for the general case. We define {\em Igusa invariants}, $J_{2i}(E)$, $i=1,2,3,5$, of the Weierstrass equation $E$ to be the ones associated to $P(x)+Q(x)^2/4$.

 Another set of invariants, {\em Igusa-Clebsch invariants}, $I_2, I_4, I_6, I_{10}$ may be associated to the Weierstrass equation $E$. The Igusa invariants are linked to the Igusa-Clebsch invariants as follows:
\begin{eqnarray}\label{eq2}
J_2 &=& I_2 / 8,\nonumber \\
J_4 &=& (4 J_2^2 - I_4) / 96,\nonumber \\
J_6 &=& (8 J_2^3 - 160 J_2 J_4 - I_6) / 576, \nonumber\\
J_8 &=& (J_2 J_6 - J_4^2) / 4,\nonumber \\
J_{10} &=& I_{10} / 4096.
\end{eqnarray}

 If $E':y'^2+S(x')y'=T(x')$ is another Weierstrass equation describing $C$, then there exist $a,b,c,d\in K$, $e\in K^{\times}$, $H[x]\in K[x]$ such that
  \begin{eqnarray}\label{eq3}x'=\frac{ax+b}{cx+d},\qquad y'=\frac{ey+H(x)}{(cx+d)^3},\qquad \textrm{ where }ad-bc\ne0.\end{eqnarray} 
In this case, we say that $E$ and $E'$ are {\em $K$-equivalent}.   Two $K$-equivalent Weierstrass equations describe isomorphic genus 2 curves over $K$. 
  
  If $\Char K\ne 2$, then one may assume that $Q(x)=0$. Further, one can see that \begin{eqnarray}\label{eq4}J_{10}(E')=J_{10}(E)e^{-20}(ad - bc)^{30}.\end{eqnarray}
   In fact, two Weierstrass equations $E$ and $E'$ describe isomorphic genus two curves over $\overline K$ if and only if there is a $\lambda\in\overline{ K}^{\times}$ such that $J_{2i}(E')=\lambda^{2i} J_{2i}(E)$.

  We now assume that $K$ is a discrete valuation field with ring of integers $\OK$, and normalized discrete valuation $\nu$. If $C$ is a genus two curve defined over $K$ by a Weierstrass equation $E:y^2+Q(x)y=P(x)$, then $E$ is said to be {\em integral} if $P(x),Q(x)\in\OK[x]$. An integral Weierstrass equation $E$ is {\em minimal} over $K$ if $\nu(J_{10}(E))\le \nu(J_{10}(E'))$ for any other integral Weierstrass equation $E'$ describing $C$ over $K$. More precisely, $\nu(J_{10}(E'))=\nu(J_{10}(E))+10 m$, where $m\ge 0$ is given by the formula in (\ref{eq4}).

  The curve $C$ is said to have {\em good reduction} over $K$ if there exists an integral Weierstrass equation $E$ describing $C$ over $K$ such that $\nu(J_{10}(E))=0$. It is said to be of {\em potential good reduction} if there exists a finite extension $K\subset K'$ such that $C$ attains good reduction over $K'$. The following criterion for potential good reduction is due to Igusa, \cite{Igusa}.
  \begin{Theorem}
  \label{aux1}
  Let $C$ be a genus two curve defined over $K$ by the Weierstrass equation $E:y^2+Q(x)y=P(x)$. The curve $C$ has potential good reduction if and only if
$J_{2i}(E)^5/J_{10}(E)^i\in \OK$ for every $i\le 5$.
  \end{Theorem}

  Assuming that $K$ is a number field with ring of integers $\OK$. A Weierstrass equation $E:y^2+Q(x)y=P(x)$ describing a genus two curve $C$ defined over $K$ is said to be {\em integral} if $P(x),Q(x)\in\OK[x]$; whereas an integral Weierstrass equation $E$ is said to be {\em minimal} at a prime $\mathfrak{p}$ of $K$ if it is minimal over $K_{\mathfrak p}$, the completion of $K$ at $\mathfrak{p}$. The equation $E$ is said to be {\em globally minimal} if it is minimal at every prime $\mathfrak{p}$ of $K$. If $K$ has class number 1, then one can always describe $C$ by a globally minimal Weierstrass equation, see \cite[Proposition 2]{Liu}.

   Similarly, $C$ has good reduction at a prime $\mathfrak{p}$ of $K$ if it has good reduction over $K_{\mathfrak{p}}$, and $C$ is said to have {\em everywhere good reduction} if $C$ has good reduction at every prime $\mathfrak{p}$ of $K$.

    Finally, one knows that if a genus two curve has good reduction at a prime $\mathfrak p$ of $K$, then its Jacobian has good reduction at $\mathfrak p$.

\section{Mestre's algorithm}
\label{sec:Mestre}
Let $K$ be a number field. If $E$ is a Weierstrass equation describing a genus 2 curve such that the Igusa invariants $J_2(E),J_4(E),J_6(E)$ and $J_{10}(E)$ are $K$-rational, then this does not mean necessarily that $C$ can be described by a Weierstrass equation over $K$. In \cite{Mestre}, Mestre was able to give a necessary and sufficient condition for $C$ to be defined by a Weierstrass equation over $K$ when the only nontrivial element in the automorphism group of $C$ is the hyperelliptic involution. For genus 2 curves with different automorphism groups, it was shown in \cite{CardonaQuer} that if the Igusa invariants of $E$ are $K$-rational, then $E$ is $\overline K$-equivalent to a Weierstrass equation defined over $K$.

In what follows, we explain the construction of Mestre and Cardona-Quer. The description of the construction is algorithmic and it is implemented in \Magma \cite{Magma} and \Sage \cite{SteinSage}. Given a quadruple $(J_2,J_4,J_6,J_{10})\in K^4$, $J_{10}\ne 0$, we want to decide whether there is a Weierstrass equation 
 defined over $K$ whose Igusa invariants are given by this quadruple. Equivalently, one may consider the following problem. If $F$ is a field with characteristic $\ne 2,3,5$, we set
\[\mathcal{M}_2(\overline F)=\{(J_2,J_4,J_6,J_{10})\in \overline{F}^4:J_{10}\ne 0\}/F^{\times}\]
where the action of $F^{\times}$ is defined by a weighted scaling of the form $$a\cdot(J_2,J_4,J_6,J_{10})=(a^2J_2,a^4J_4,a^6J_6,a^{10}J_{10}).$$ We write $(J_2',J_4',J_6',J_{10}')\sim_F(J_2,J_4,J_6,J_{10})$ if $(J_2',J_4',J_6',J_{10}')=a\cdot(J_2,J_4,J_6,J_{10})$ for some $a\in F^{\times}$. A {\em moduli point} $J\in\mathcal{M}_2(\overline F)$ is {\em defined over $F$} if it is the equivalence class of a quadruple $(J_2,J_4,J_6,J_{10})$ with all $J_{2i}\in F$, $i=1,2,3,5$.

Let $C$ be a curve of genus 2 defined over $\overline{K}$ by a Weierstrass equation $E$. We say that $F\subset\overline{K}$ is a {\em field of moduli} for $C$ if the equivalence class of $(J_2(E),J_4(E),J_6(E),J_{10}(E))$ in $\mathcal{M}_2(\overline K)$ is defined over $F$. We say that $F\subset\overline{K}$ is a {\em field of definition} for $C$ if there exists a Weierstrass equation $E'$ defined over $F$ which is $\overline K$-equivalent to $E$.  The problem now boils down to: given that $K$ is a field of moduli for $C$, is $K$ a field of definition for $C$?

 Let $C$ be a genus 2 curve such that $K$ is a field of moduli for $C$. One may construct homogeneous
ternary forms $Q_J,T_J\in K[u,v,w]$ of degrees 2 and 3 respectively, see \cite[\S 1.5]{Mestre} for the explicit description and properties of $Q_J$ and $T_J$. The conic $M_J$ defined by $Q_J=0$ has a $K$-rational point if $C$ has a nontrivial automorphism other than the hyperelliptic involution, see \cite{CardonaQuer}. On the other hand, the conic defined by $M_J:Q_J=0$ may not have a $K$-rational point if the only nontrivial automorphism of $C$ is the hyperelliptic involution. In the latter case, if $M_J(K)=\emptyset$, then there exists a quadratic extension $K\subset K'$ such that $M_J(K')\ne\emptyset$, see \cite{Mestre}. The existence of a rational point over $L\in\{K,K'\}$ gives rise to a parametrization $(u,v,w)=(u(x),v(x),w(x))$ where $u,v,w\in L[x]$ are polynomials of degree 2. Now we get the Weierstrass equation $y^2=T_J(u(x),v(x),w(x))$ describing $C$ where $T_J(u(x),v(x),w(x))\in L[x]$ is a degree 6 polynomial. The following theorem summarizes the above discussion.
\begin{Theorem}
\label{Aux2}
Let $E$ be a Weierstrass equation describing a genus 2 curve $C$. Assume that $J_2(E),J_4(E),J_6(E),J_{10}(E)\in K$.  Let $J$ be the class of the quadruple of Igusa invariants $(J_2(E),J_4(E),\\ J_6(E),J_{10}(E))$ of $E$ in $\mathcal{M}_2(\overline K)$. 
\begin{itemize}
\item[(1)] (Mestre, \cite{Mestre}) If $\Aut(C)\isom \Z_2$, then
$K$ is a field of definition for $C$ if and only if $M_J(K)\ne \emptyset$.
\item[(2)] (Cardona-Quer, \cite{CardonaQuer}) If $\Aut(C)\not\isom \Z_2$, then $K$ is always a field of definition for $C$.
\end{itemize}
\end{Theorem}
One remarks that when $\Aut(C)\isom\Z_2$, if $M_J(K)=\emptyset$, there are many quadratic extensions $K\subset K'$ over which $M_J(K')\ne\emptyset$. We set two of the coordinates $u,v,w$ to be fixed and solve for the third coordinate. This gives a rational point on the conic over a quadratic extension $K\subset K'$.

 One disadvantage about Mestre's algorithm is that the coefficients of the sextic $T_J(u(x),v(x),w(x))$ are large numbers. This is why we will spend the next section discussing minimization of genus two curves.

\section{Minimization of genus two curves} 

In this section, $K$ is a complete discrete valuation field with ring of integers $\OK$, maximal ideal $\mathfrak{m}$, residue field $k=\OK/\mathfrak m$, a uniformiser $\pi$, and normalized discrete valuation $\nu$. We assume that $\Char k>5$ throughout this section. Moreover, we will be mainly concerned with genus two curves described by Weierstrass equations $E$ defined over $K$ with potential good reduction. In other words, $J_{2i}^5(E)/J_{10}^i(E)\in \OK$ for every $i\le 5$, see Theorem \ref{aux1}.

Let $(j_2,j_4,j_6,j_{10})\in K^4$ where $j_{10}\ne 0$. Let $C$ be a genus 2 curve defined by a Weierstrass equation $E$ over $K$ with $(J_2(E),J_4(E),J_6(E),J_{10}(E))=\lambda\cdot(j_2,j_4,j_6,j_{10})$ for some $\lambda\in K^{\times}.$ Assume, moreover, that $C$ has potential good reduction over $K$. In what follows, we discuss the possibilities for $\nu(J_{10}(E))$ if $E$ is a minimal Weierstrass equation over $K$.
\begin{Proposition}
  \label{prop1}
  Let $E_1$ and $E_2$ be two minimal Weierstrass equations over $K$ describing two genus 2 curves with potential good reduction over $K$. Assume moreover that $E_1$ and $E_2$ are $\overline{K}$-equivalent.
 Then $$|\nu(J_{10}(E_1))-\nu(J_{10}(E_2))|\in\{0,5,10,15,20\}.$$
\end{Proposition}
\begin{Proof}
We set $\overline{\nu}$ to be an extension of $\nu$ to $\overline{K}$.
   One knows that $(J_2(E_2),J_4(E_2),J_6(E_2),J_{10}(E_2))=(\lambda^2 J_2(E_1),\lambda^4 J_4(E_1),\lambda^6 J_6(E_1),\lambda^{10} J_{10}(E_1))$ for some $\lambda\in\overline{K}^{\times}$. The fact that both $E_1$, $E_2$ are minimal, hence integral, and $J_2(E_2)=\lambda^2J_2(E_1)\in \OK$, imply that $\lambda^2\in K^{\times}$. Since both curves have potential good reduction over $K$ and $\Char k>5$, this implies that $\nu(J_{10}(E_1)),\nu(J_{10}(E_2))\le 20$, see \cite[Remarque 13]{Liu}. Now, since $\nu(J_{10}(E_2))=10\overline{\nu}(\lambda)+\nu(J_{10}(E_1))\le 20$, it follows that
 $10\overline{\nu}(\lambda)\in\{0,5,10,15,20\}$.
\end{Proof}

An integral Weierstrass equation $E:y^2+Q(x)y=P(x)$ over $K$ will be called {\em reduced} over $K$ if $4P(x)+Q(x)^2\not\equiv 0$ in $k[x]$, otherwise, it is called {\em non-reduced} over $K$.
\begin{Remark}
\label{rem:reduced}
It is easy to see that if $E$ is an integral Weierstrass equation describing a genus 2 curve $C$ over $K$, and $E$ is non-reduced over $K$, then the transformation $x\mapsto x $, $y\mapsto \pi^{1/2}y$ yields an integral Weierstrass equation $E'$ over $K$ that describes a genus 2 curve isomorphic to $C$ over $K$, and $\nu(J_{10}(E'))=\nu(J_{10}(E))-10$.

We also remark that if $E$ is a minimal Weierstrass equation over $K$ with $\nu(J_{10}(E))\le9$, then $E$ is reduced over $K$, see \cite[Proposition 7]{Liu}.   
\end{Remark}
\begin{Remark}
\label{rem:equivalenttoreduced}
Let $E:y^2+Q(x)y=P(x)$ be integral over $K$ and $4P(x)+Q(x)^2=\sum_{i=0}^6a_ix^i$. If $$\nu(a_0)\ge 7,\; \nu(a_1)\ge 5,\; \nu(a_2)\ge 3,\; \nu(a_3)\ge 1,$$ then $E$ is $K$-equivalent to a non-reduced Weierstrass equation $E'$ via the transformation  
$x\mapsto \pi^2 x $, $y\mapsto \pi^{3}y$. Moreover, $\nu(J_{10}(E'))=\nu(J_{10}(E))$.
\end{Remark}
We recall that if $E$ is a minimal Weierstrass equation over $K$ describing a genus 2 curve $C$ with $J_2(E),J_4(E),J_6(E)\in\OK$ and $J_{10}(E)\in\mathcal{O}_K^{\times}$, then the curve $C$ has potential good reduction since $J_{2i}^5(E)/J_{10}^i(E)\in\OK$ for every $i$, see Theorem \ref{aux1}. The fact that $C$ has potential good reduction implies the existence of a finite extension $F$ of $K$ over which $C$ attains good reduction. Moreover, one may assume that $F$ is minimal, i.e., $F$ is contained in any other such field, see \cite[\S3]{Liu3}. Now we state the main result of this section. 

\begin{Theorem}
\label{thm:local} Let $J_2,J_4,J_6,J_{10}$ be such that $J_i\in \OK$, $i=2,4,6$, and $J_{10}\in\mathcal{O}_K^{\times}$. Let $E$ be a minimal Weierstrass equation over $K$ describing a genus 2 curve $C$ with $(J_2(E),J_4(E),J_6(E),J_{10}(E))\sim_{\overline{K}}(J_2,J_4,J_6,J_{10})$. Let $F$ be the minimal field extension over which $C$ attains good reduction and suppose that $[F:K]>1$.   

If $E$ is not $F$-equivalent to a non-reduced Weierstrass equation $E'$ over $F$ with $\nu_F(J_{10}(E'))=\nu_F(J_{10}(E))$, then $\nu(J_{10}(E))\in\{5,10,15\}.$ Moreover, one has
\begin{itemize}
\item[a)]  If $\nu(J_{10}(E))=5$, then $[F:K]=6$.
\item[b)] If $\nu(J_{10}(E))=10$, then $[F:K]=3$.
\item[c)] If $\nu(J_{10}(E))=15$, then $[F:K]= 2$.
\end{itemize}
\end{Theorem} 

\begin{Proof} 
 One sees that \((J_2(E),J_4(E),J_6(E),J_{10}(E))= (\lambda^2J_2,\lambda^4J_4,\lambda^6J_6,\lambda^{10}J_{10})\) for some $\lambda\in \overline{K}^{\times}$. Now since $E$ is minimal hence integral, it follows that $\lambda^2 J_2\in\OK$. Therefore, $\lambda^2\in K$. Since $\nu(J_{10}(E))=5\nu(\lambda^2)+\nu(J_{10})=5\nu(\lambda^2)\le 20$, \cite[Remarque 13]{Liu}, it follows that $\nu(J_{10}(E))\in\{0,5,10,15,20\}$. We notice that $\nu(J_{10}(E))\ne0$ as $[F:K]>1$.

   Since $C$ has potential good reduction over $K$, and $\Char k>5$, it follows that the minimal field $F$ over which $C$ attains good reduction is a totally and tamely ramified extension of $K$, hence cyclic, with $n_F:=[F:K]\in\{2,3,4,5,6,8,10\}$, see \cite[Corollaire 4.1]{Liu2}. We set $\pi_F$ to be a uniformiser in the ring of integers $\mathcal{O}_F$ of $F$, and we assume without loss of generality that $\pi_F^{n_F}$ is either $\pi$ or $-\pi$. We set $\nu_F$ to be the corresponding discrete valuation on $F$.
   
   One knows that since $E:y^2+Q(x)=P(x)$ is not minimal over $F$, this yields that either $E$ is non-reduced over $K$ or $4P(x)+Q(x)^2=\sum_{i=0}^6 a_ix^i\in\OK[x]$ has a root $x_0$ of multiplicity at least $4$ in $k$, see for example \cite[Lemma 5.3]{MSt}. 
   
     Since $E$ is reduced by assumption, then $4P(x)+Q(x)^2$ has a root $x_0$ of multiplicity at least 4.  We use the translation map $x-t$ where $t$ is a lift of $x_0$ in $\OK$. This map does not change $\nu(J_{10}(E))$, and it allows us to assume that $\nu(a_i)\ge1$ when $i=0,1,2,3$, as now $0$ is the root of $4P(x)+Q(x)^2$ in $k[x]$ that has multiplicity at least $4$. Now we may apply the transformation $x\mapsto \pi_{F}^r x$, $y\mapsto \pi_{F}^s y$, for some $r,s\ge1$, to obtain the minimal Weierstrass equation $E'$ over $F$ with $\nu_{F}(J_{10}(E'))=0$, see for example \cite[Lemma 11.2]{MSt}. Since $E$ and $E'$ are $F$-equivalent Weierstrass equations, the identity (\ref{eq4}) implies that $n_F\nu(J_{10}(E))=20s-30r$. In the following table, we run through the possibilities for $r$ and $s$ when $n_F\in\{2,3,4,5,6,8,10\}$. The entries of the table contain the solutions $(r,s)$ for the linear equation $n_F\nu(J_{10}(E))=20s-30r$ where $t\ge 0$ is an integer. 
      
      \begin{center}
      {\tiny
        \begin{tabular}{|c|c|c|c|c|c|c|c|}
        \hline
       \diagbox{ $\nu(J_{10}(E))$}{$n_F$}& 2 &3& 4& 5 & 6 & 8 & 10  \\
          \hline
          5& $(1+2t,2+3t)$& -&$(2+2t,4+3t)$&-& $(1+2t,3+3t)$& $(2+2t,5+3t)$& $(1+2t,4+3t)$ \\
          \hline
          10&$(2+2t,4+3t)$&$(1+2t,3+3t)$& $(2+2t,5+3t)$&$(1+2t,4+3t)$&$(2+2t,6+3t)$& $(2+2t,7+3t)$& $(2+2t,8+3t)$ \\
          \hline
          15&$(1+2t,3+3t)$& -&$(2+2t,6+3t)$&-&$(1+2t,6+3t)$& $(2+2t,9+3t)$& $(1+2t,9+3t)$ \\
          \hline
          20&$(2+2t,5+3t)$&$(2+2t,6+3t)$& $(2+2t,7+3t)$& $(2+2t,8+3t)$& $(2+2t,9+3t)$& $(2+2t,11+3t)$& $(2+2t,13+3t)$\\
          \hline
        \end{tabular}}
      \end{center}
  
  One remarks that $E'$ is given by $y^2+Q_1(x)y=P_1(x)$, thus 
  \begin{eqnarray}\label{eqq}4P_1(x)+Q_1(x)^2=\sum_{i=0}^6 \pi_F^{ir-2s}a_ix^i\in\mathcal{O}_F[x].\end{eqnarray} Since $E'$ describes a curve with good reduction over $F$, it follows that $4P_1(x)+Q_1(x)^2$ is separable of degree at least $5$ in $k_F[x]$. Now, if $5r-2s>0$, then this means that the latter polynomial is of degree at most $4$ in $k_F[x]$. Also, we note that if $5r-2s<0$, then the integrality of $E'$ implies $\nu(a_5)>\lfloor(2s-5r)/n_F\rfloor$; similarly, if $6r-2s<0$, then $\nu(a_6)>\lfloor(2s-6r)/n_F\rfloor$. Thus,  if $6r-2s<0$, then $4P_1(x)+Q_1(x)^2$ is again of degree at most $4$ in $k_F[x]$. This leaves us with the following possibilities:
  
  \begin{center}
      {\tiny
        \begin{tabular}{|c|c|c|c|c|c|c|c|}
        \hline
       \diagbox{ $\nu(J_{10}(E))$}{$n_F$}& 2 &3& 4& 5 & 6 & 8 & 10  \\
          \hline
          5& -& -&-&-& $(1,3)$& $(2,5)$& - \\
          \hline
          10&-&$(1,3)$& $(2,5)$&-&$(2,6)$& $(2+2t,7+3t)$, $t=1$& $(2+2t,8+3t)$, $t=1$ \\
          \hline
          15&$(1,3)$& -&-&-&$(1,6)$& $(2+2t,9+3t)$, $t=0,1,2$& $(1+2t,9+3t)$, $t=2,3$ \\
          \hline
          20&$(2,5)$&$(2,6)$& $(2+2t,7+3t)$, $t=0,1$& -& $(2+2t,9+3t)$, $t=0,1,2$& $(2+2t,11+3t)$, $t=1,3$& $(2+2t,13+3t)$, $t=4$\\
          \hline
        \end{tabular}}
      \end{center}
      
      We now consider the case $\nu(J_{10}(E))=5$. 
If $n_F=8$, the pair $(r,s)$ takes the value $(2,5)$.  In the latter case, we have     $$\nu_F(a_0)\ge16,\ \nu_F(a_1)\ge8,\ \nu_F(a_2)\ge8,\ \nu_F(a_3)\ge 8,\ \nu_F(a_4)\ge8.$$ Hence Remark \ref{rem:equivalenttoreduced} implies that $E$ is $F$-equivalent to a non-reduced Weierstrass equation $E'$ over $F$ with $\nu_F(J_{10}(E'))=\nu_F(J_{10}(E))$.

We consider the case $\nu(J_{10}(E))=10$. If $n_F=4$,  the pair $(r,s)$ takes the value $(2,5)$. The integrality of $E$ over $K$ implies that $$\nu_F(a_0)\ge12,\ \nu_F(a_1)\ge8,\ \nu_F(a_2)\ge8,\ \nu_F(a_3)\ge 4,\ \nu_F(a_4)\ge4.$$ We conclude by Remark \ref{rem:equivalenttoreduced} that this pair does not allow $E$ to satisfy the hypothesis of the theorem. For $n_F=6$ and the pair $(r,s)=(2,6)$, we conclude similarly as then    $$\nu_F(a_0)\ge12,\ \nu_F(a_1)\ge12,\ \nu_F(a_2)\ge12,\ \nu_F(a_3)\ge 6,\ \nu_F(a_4)\ge6,\ \nu_F(a_5)\ge6.$$ In a similar fashion, one has for $n_F=8$ and the pair $(r,s)=(4,10)$ the following valuations  $$\nu_F(a_0)\ge24,\ \nu_F(a_1)\ge16,\ \nu_F(a_2)\ge16,\ \nu_F(a_3)\ge 8,\ \nu_F(a_4)\ge8,$$ whereas for $n_F=10$ and the pair $(r,s)=(4,11)$,  one has $$\nu_F(a_0)\ge30,\ \nu_F(a_1)\ge20,\ \nu_F(a_2)\ge20,\ \nu_F(a_3)\ge 10,\ \nu_F(a_4)\ge10.$$ 

A similar argument shows that when $\nu(J_{10}(E))=15$, then the only admissible pair is $(1,3)$ when $n_F=2$, otherwise for $n_F\ge 6$ and a corresponding pair from the latter table will force $E$ to be $F$-equivalent to a non-reduced Weierstrass equation $E'$ over $F$ with $\nu_F(J_{10}(E'))=\nu_F(J_{10}(E))$. 
   As for $\nu(J_{10}(E))=20$, any choice for $n_F$ and a corresponding pair of $(r,s)$ will again yield that $E$ is  
      $F$-equivalent to a non-reduced Weierstrass equation $E'$ over $F$ with $\nu_F(J_{10}(E'))=\nu_F(J_{10}(E))$.   
\end{Proof}

One notices that the converse of Theorem \ref{thm:local} does not necessarily hold. For example the Weierstrass equation $E: y^2=\pi x^6+\pi $, where $\Char k>5 $, is minimal. The equation $E$ describes a genus $2$ curve that has potential good reduction. In addition, one may check easily that $\nu(J_{10}(E))=10$. However, $E$ is clearly non-reduced.

In the proof of Theorem \ref{thm:local}, we mentioned the minimization transformations for $E$ over the minimal field $F$ to produce an integral equation $E'$ over $F$ with $\nu_{F}(J_{10}(E'))=0$. We list these transformations in the following proposition.

\begin{Proposition}
\label{prop:minimization}
Let $J_2,J_4,J_6,J_{10}$ be such that $J_i\in \OK$, $i=2,4,6$, and $J_{10}\in\mathcal{O}_K^{\times}$. Let $E:y^2+Q(x)y=P(x)$ be a minimal Weierstrass equation over $K$ describing a genus 2 curve $C$ with $(J_2(E),J_4(E),J_6(E),J_{10}(E))\sim_{\overline{K}}(J_2,J_4,J_6,J_{10})$. Let $F$ be the minimal field extension over which $C$ attains good reduction and assume that $[F:K]>1$.   

If $E$ is not $F$-equivalent to a non-reduced Weierstrass equation $E'$ over $F$ with $\nu_F(J_{10}(E'))=\nu_F(J_{10}(E))$, then $E$ is $F$-equivalent to an integral Weierstrass equation $E''$ over $F$ with $\nu_F(J_{10}(E''))=0$. Moreover, $4P(x)+Q(x)^2$ has a root of multiplicity $6$ in $k[x].$
\end{Proposition}
\begin{Proof}
Set $n_F=30/\nu(J_{10}(E))$ and $\pi_F=(\pm\pi)^{1/n_F}$. The polynomial $4P(x)+Q(x)^2$ has a root of multiplicity at least $4$ in $k[x]$. After a translation one may assume that this root is $0$. Moreover, Theorem \ref{thm:local} together with equality (\ref{eqq}) imply that this root is of multiplicity $6$. 
Now, one may use the transformation $(x,y)\mapsto(\pi_{F}x, \pi_{F}^3y)$ to obtain an integral Weierstrass equation $E'$ over $F$ with $\nu_F(J_{10}(E'))=0$. 
\end{Proof}
\begin{Remark}
\label{rem:Fequivalentnonreduced}
In Theorem \ref{thm:local}, one may relax the condition that $E$ is not $F$-equivalent to a non-reduced Weierstrass equation $E'$ over $F$ with $\nu_F(J_{10}(E'))=\nu_F(J_{10}(E))$ and replace it with $E$ is not $K$-equivalent to a non-reduced Weierstrass equation $E'$ over $K$ with $\nu(J_{10}(E'))=\nu(J_{10}(E))$. This happens because in the three cases of the theorem, one has $\nu_F(J_{10}(E))=30$,  hence if $E$ is $F$-equivalent to a non-reduced Weierstrass equation $E'$ over $F$ with $\nu_F(J_{10}(E'))=\nu_F(J_{10}(E))$, then one may replace $E$ with an integral equation $E_1$ over $F$ with $\nu_F(J_{10}(E_1))=20$. But since the curve $C$ attains good reduction over $F$, then one may find $r$ and $s$ such that $20s-30r=20$, i.e., $(r,s)=(2+2t,4+3t)$ but as in the proof of the theorem, this implies that $C$ does not have good reduction over $F$.
\end{Remark}

Let $C$ be a genus $2$ curve over $K$. Let $\mathcal{C}$ be a stable model of $C$ over the strict henselisation of $\OK$ in the sense of \cite[D\'{e}finition 1]{Liu2}. Let $f:\mathcal C\to\mathcal Z:=\mathcal C/\langle\sigma\rangle$ where $\sigma$ is the extension of the hyperelliptic involution of $C$ to $\mathcal C$. We set $\omega$ to be the point in the generic fiber of $\mathcal Z$ corresponding to $x=\infty$, and $\overline{\omega}$ its specialization in the special fiber of $\mathcal Z$. Then $\overline{\omega}$ is said to be ramified if $f$ is ramified at $\overline{\omega}$; otherwise $\overline{\omega}$ is said to be non-ramified. 

\begin{Corollary}
\label{cor:irreducible}
Let $J_2,J_4,J_6,J_{10}$ be such that $J_i\in \OK$, $i=2,4,6$, and $J_{10}\in\mathcal{O}_K^{\times}$. Let $E$ be a minimal Weierstrass equation over $K$ describing a genus 2 curve $C$ with $(J_2(E),J_4(E),J_6(E),J_{10}(E))\sim_{\overline{K}}(J_2,J_4,J_6,J_{10})$. 
Let $F$ be the minimal field extension over which $C$ attains good reduction and suppose that $[F:K]>1$.   

If $E$ is not $F$-equivalent to a non-reduced Weierstrass equation $E'$ over $F$ with $\nu_F(J_{10}(E'))=\nu_F(J_{10}(E))$, then $\overline{\omega}$ is non-ramified.
\end{Corollary}
\begin{Proof}
 The fact that $n_F$ is either $6,3$, or $2$ when $\nu(J_{10}(E))=5,10$, or $15$, respectively, see Theorem \ref{thm:local}, implies that $\overline{\omega}$ is non-ramified, see \cite[Th\'eor\`{e}me 1]{Liu3}.
\end{Proof}

\section{Globally minimal Weierstrass equations} 

In this section we show that the existence of a minimal Weierstrass equation describing a genus 2 curve over a number field $K$ whose Igusa invariant $J_{10}$ satisfies certain properties implies the existence of a genus two curve with everywhere good reduction over a quadratic extension of $K$.

We assume throughout this section that $K$ is a number field with ring of integers $\OK$. If $\mathfrak p$ is a prime of $K$, then $K_{\mathfrak p}$ is the completion of $K$ at $\mathfrak p$, with ring of integers $\mathcal{O}_{K_{\mathfrak p}}$, maximal ideal $\mathfrak{m_{\mathfrak p}}$, residue field $k_{\mathfrak p}$, a uniformiser $\pi_{\mathfrak p}$ in $\OK$, and normalized discrete valuation $\nu_{\mathfrak p}$.

\begin{Theorem}
\label{thm1}
Let $K$ be a number field with class number $1$. Let $(J_2,J_4,J_6,J_{10})$ be such that $J_i\in \OK$, $i=2,4,6$, and $J_{10}\in\mathcal{O}_K^{\times}$. Let $E$ be a globally minimal Weierstrass equation over $K$ describing a genus 2 curve $C$ with $(J_2(E),J_4(E),J_6(E),J_{10}(E))\sim_{\overline K}(J_2,J_4,J_6,J_{10})$.
 Assume that for every bad prime $\mathfrak p$ of $C$ in $K$, one has $\gcd(\Char k_{\mathfrak{p}},30)=1$, and moreover, 
one of the following conditions is satisfied:
 \begin{itemize}
 \item[i)] either $\nu_{\mathfrak p}(J_{10}(E))=10$ and $E$ is non-reduced over $K_{\mathfrak{p}}$,
 \item[ii)] or $\nu_{\mathfrak{p}}(J_{10}(E))=15$, and $E$ is not $K_{\mathfrak{p}}$-equivalent to a non-reduced Weierstrass equation $E'$  over $K_{\mathfrak{p}}$ with $\nu_{\mathfrak p}(J_{10}(E'))=15$.
 \end{itemize} Then $C$ has everywhere good reduction over $L:=K\left((\pm\pi_{{\mathfrak p}_1}\ldots\pi_{{\mathfrak p}_k})^{1/2}\right)$ where $\mathfrak{p}_i$'s are the primes such that $\nu_{\mathfrak{p}_i}(J_{10}(E))=15$. Moreover, $C$ is described by a globally minimal Weierstrass equation over $L$. 
\end{Theorem}

\begin{Proof}
In view of Remark \ref{rem:reduced},  if $E:y^2+Q(x)y=P(x)$ is non-reduced over $K_{\mathfrak{p}}$ and $\nu_{\mathfrak p}(J_{10}(E))=10$, then we apply the transformation $x\mapsto x$ and $y\mapsto \pi_{\mathfrak{p}}^{1/2} y$.  This yields an integral Weierstrass equation $E'$ over $K_{\mathfrak{p}}$ with $\nu_{\mathfrak p}(J_{10}(E'))=0$.

 So we may (and will) assume that every bad prime of $C$ is of type ii).  Now that $\gcd(\Char k_{\mathfrak{p}},30)=1$, Proposition \ref{prop:minimization} asserts that for every such prime $\mathfrak{p}$ one may use a translation $x\mapsto x+x_{\mathfrak p}$ where $x_{\mathfrak p}$ is an integral lift of the multiplicity-$6$ root of $4P(x)+Q(x)^2$ in $k_{\mathfrak p}$. A Chinese Remainder Theorem argument shows that there exists $x_0\in\OK$ such that $x_0$ is a lift of $x_{\mathfrak p}$ for every bad prime $\mathfrak p$. So using the transformation $x\mapsto x+x_0$, which neither changes $\nu_{\mathfrak p}(J_{10}(E))$ for any $\mathfrak p$ nor the integrality of $E$ over $K$, will allow us to use the minimization transformation $x\mapsto(\pm\pi_{\mathfrak p})^{1/2}x$, $y\mapsto(\pm\pi_{\mathfrak p})^{3/2}y$ for every bad prime $\mathfrak p$ of $C$ to obtain an integral Weierstrass equation $E_{\mathfrak p}$ over $F_{\mathfrak p}:=K_{\mathfrak p}((\pm\pi_{\mathfrak p})^{1/2})$ with $\nu_{F_{\mathfrak p}}(J_{10}(E_{\mathfrak p}))=0$, see Proposition \ref{prop:minimization}.  Therefore, the choice of $\pi_\mathfrak p$ in $\OK$ shows that using the transformation $$x\mapsto(\pm\pi_{\mathfrak p_1}\pi_{\mathfrak p_2}\ldots \pi_{\mathfrak p_k})^{1/2}x,\qquad  y\mapsto(\pm\pi_{\mathfrak p_1}\pi_{\mathfrak p_2}\ldots \pi_{\mathfrak p_k})^{3/2}y$$ will produce an integral Weierstrass equation $E'$ over $L=K\left((\pm\pi_{{\mathfrak p}_1}\ldots\pi_{{\mathfrak p}_k})^{1/2}\right)$ where ${\mathfrak p}_1,\ldots,{\mathfrak p}_k$ are the bad primes of $C$ with $\nu_{\mathfrak{p}}(J_{10}(E))=15$, and $E'$ describes a genus two curve over $L$ with everywhere good reduction.
\end{Proof}
We recall that any Weierstrass equation $E$ satisfying the assumptions of Theorem \ref{thm1} must satisfy that $\overline{\omega}$ is non-ramified, see Corollary \ref{cor:irreducible}.
\begin{Remark}
\label{rem:classnumber}
The assumption that $K$ has class number $1$ in Theorem \ref{thm1} above ensures the existence of a globally minimal Weierstrass equation describing $C$ over $K$, see \cite[Remarque 6]{Liu}. However, the existence of a globally minimal Weierstrass equation for $C$ over $K\left((\pm\pi_{{\mathfrak p}_1}\ldots\pi_{{\mathfrak p}_k})^{1/2}\right)$ is not guaranteed.  
\end{Remark}

\begin{Remark}
In \cite{Kida}, it was shown that there is no elliptic curve defined over the field of rational numbers that attains good reduction at every finite place under quadratic base change. Theorem \ref{thm1} shows that a genus two curve defined over $\Q$ may attain good reduction everywhere over a quadratic field extension of $\Q$. In fact, all genus 2 curves constructed in the coming sections are defined over $\Q$ and they assume everywhere good reduction over a quadratic field extension.
\end{Remark}

\section{Genus two curves with everywhere good reduction over quadratic fields}

In this section we describe the algorithm to compute genus 2 curves with everywhere good reduction over a quadratic extension of $\Q$. Then we list some of the Weierstrass equations $E$ provided by the algorithm. It is worth mentioning that in this section we will be producing such equations with $J_2(E)\ne 0$, whereas we will leave the case of $J_2(E)=0$ for \S \ref{sec:infinite}.

Let $J=(J_2,J_4,J_6,\pm 1)$, where $J_i\in\Z$, be such that either $J_2=0$; or $J_2\ne 0$ and $M_J(\Q)\ne\emptyset$. We recall that according to \S \ref{sec:Mestre}, the latter condition on $J_2$ implies the existence of a genus two curve $C$ defined by a globally minimal Weierstrass equation $E:y^2+Q(x)y=P(x)$ over $\Q$ with $(J_2(E),J_4(E),J_6(E),J_{10}(E))\sim_{\overline{\Q}}(J_2,J_4,J_6,\pm 1)$. If, moreover, for every prime $p\mid J_{10}(E)$, one either has $\nu_p(J_{10}(E))=10$ and $E$ is non-reduced over $\Q_p$; or $\gcd(J_{10}(E),30)=1$ and $\nu_p(J_{10}(E))=15$, then the quadruple $J$ will be said to {\em satisfy condition (*)}.

\begin{Algorithm}
\label{alg}
\end{Algorithm}
\textbf{Input:} A quadruple $J=(J_2,J_4,J_6,\pm 1)\in \Z^4$ that satisfies condition (*).

\textbf{Output:}
\begin{itemize}
\item[i)] A quadratic field extension $K$ of $\Q$.
\item[ii)] A genus two curve $C$ with everywhere good reduction over $K$ described by an integral Weierstrass equation $E$ over $K$ and $(J_2(E),J_4(E),J_6(E),J_{10}(E))\sim_{\overline{K}}(J_2,J_4,J_6,\pm 1)$. Moreover, $E$ is globally minimal except possibly at the primes of $K$ lying above the prime $2$.
\end{itemize}
The fact that $J$ satisfies condition (*) implies the existence of a globally minimal Weierstrass equation $E:y^2+Q(x)y=P(x)$ with
$$
(J_2(E),J_4(E),J_6(E),J_{10}(E))=(\lambda^2J_2,\lambda^4J_4,\lambda^6J_6,\pm\lambda^{10})\in\Q^4
$$
for some $\lambda \in \overline{\Q}$ with $[\Q(\lambda):\Q]\le2$. Moreover,  for every prime $p\mid J_{10}(E)$, one either has $\nu_p(J_{10}(E))=10$ and $E$ is non-reduced over $\Q_p$; or $\gcd(J_{10}(E),30)=1$ and $\nu_p(J_{10}(E))=15$. We pick such a Weierstrass equation $E$. 

The steps of the algorithm are as follows:
\begin{itemize}
\item[1)] For every prime $p$ with $\nu_p(J_{10}(E))=10$ and $E$ is non-reduced over $\Q_p$, we replace $E$ with the Weierstrass equation $ y^2=(4P(x)+Q(x)^2)/p$. We remark that the invariant $J_{10}$ associated to the latter Weierstrass equation has the same valuation as $J_{10}(E)$ at any  prime different from $ p$, whereas its $p$-valuation is zero.
     \item[2)] For every prime $p$ with $\nu_p(J_{10}(E))=15$, the polynomial $4P(x)+Q(x)^2$ has a root $x_p$ of multiplicity $6$ in $\mathbb{F}_p[x]$. Let $x_0\in\Z$ be a lift of $x_p$ for every such $p$. Now, we replace the Weierstrass equation $E$ by the Weierstrass equation obtained from $E$ using the transformation $x\mapsto x+x_0$. We recall that the latter transformation does not change $J_{10}(E)$.
    \item[3)] If $p_1,\ldots,p_n$ is the list of primes in 2), then use the transformation $$x\mapsto (\pm p_1\ldots p_n)^{1/2}x,\qquad y\mapsto (\pm p_1\ldots p_n)^{3/2}y.$$ The latter transformation gives rise to an integral Weierstrass equation over $K=\Q(\sqrt{\pm p_1\cdots p_n})$, where $p_1,\cdots,p_n$ are the prime divisors of $J_{10}(E)$ with $\nu_p(J_{10}(E))=15$. Moreover, the latter Weierstrass equation has everywhere good reduction over $K$, see Theorem \ref{thm1}.
    \end{itemize}
    One remarks that the fact that the input of Algorithm \ref{alg} satisfies condition (*) implies that the associated Weierstrass equation $E$ satisfies the hypothesis of Theorem \ref{thm1}. In particular, the algorithm must terminate with the aforementioned output as explained in the proof of Theorem \ref{thm1}.

\begin{Example} To illustrate the algorithm, we consider the quadruple $J=(J_2,J_4,J_6,J_{10}) = (1,1,1,1)$.
Using \Sage one knows that $M_J(\Q)\ne\emptyset$ and a hyperelliptic curve $C$ corresponding to this moduli point is given. Meanwhile, one may use \Magma to obtain a minimal Weierstrass equation $E$ describing $C$ over $\mathbb Q$:
{\footnotesize\begin{eqnarray*}
y^2 &+& (-x^3 + 2x^2 - 3x - 1)\;y = 319942251280530334443517978446668563494310791\;x^6  \\&-&
1348554939815014300021284658314762636324173354\;x^5 +2368396716413597071671039091048166857266661725\;x^4\\&
-&2218395544527845496581055490986877026155157575\;x^3 +1168815723018377789171495695634957245626408185\;x^2\\
& -&328436816477826084131551647893678647088163154\;x +38454430231425908654742571540402070766182706.
 \end{eqnarray*}}
The discriminant $J_{10}(E)=5^{10} \cdot 17^{15} \cdot 193^{15} \cdot 1193^{15}$. In fact, $E$ is non-reduced at $5$, so using the minimization transformations in the algorithm described above,
we obtain a genus $2$ curve with everywhere good reduction over the quadratic field
$\mathbb Q(\sqrt{ 17\cdot 193 \cdot 1193})$. A minimal Weierstrass equation $E^+_{\min}$ describing the curve is given  as follows: We set $a = \sqrt{17\cdot 193 \cdot 1193}$, and remark that $D=17\cdot 193 \cdot 1193\equiv1$ mod $4$.
{\footnotesize\begin{align*}
E^+_{\min}:\quad y^2 + (-x^3 + 2x^2 - 3x - 1)\;y &=63988450256106066888703595689333712698862158\;  x^6 \\&+ (1+14905246457204355953359713280206751857219456 \cdot a)\; x^5\\ &+5662551478017146195184525794000839333887206015356\;  x^4\\&+
\frac{1}{2}(5+586228886401597062058919441470649503237809367473\cdot a)\; x^3\\ &+
33406538328776427426641190939693462447792426559697895\;  x^2\\& +
\frac{1}{2}(-3+1037547005393506172459895771934774485127403691693305\cdot a)\; x\\ &+
13138915763144011510044623243960742638241073415065183915
 \end{align*}}
 In a similar fashion, the Weierstrass equation $E$ yields a genus $2$ curve with everywhere good reduction over the quadratic field
$\mathbb Q(\sqrt{- 17\cdot 193 \cdot 1193})$. Setting $a=\sqrt{-17\cdot 193\cdot 1193}$, a minimal Weierstrass equation $E^-_{\min}$ describing the curve is given  as follows:
{\footnotesize\begin{align*}
E^-_{\min}:\quad y^2 + (-x^3 + 2x^2 - 3x - a)\;y &=63988450256106066888703595689333712698862158 \;x^6\\&+ (1- 14905246457204355953359713280206751857219456\cdot a)\; x^5 \\& - 5662551478017146195184525794000839333887206015361\; x^4\\ &+(3 +
 293114443200798531029459720735324751618904683736\cdot a)\;   x^3\\&+(33406538328776427426641190939693462447792426559697894  +a)\; x^2\\& -
 518773502696753086229947885967387242563701845846654\cdot a \; x \\  &-13138915763144011510044623243960742638241073415064205357.
 \end{align*}}
\end{Example}
In the Appendix, we list more globally minimal Weierstrass equations $E:y^2+Q(x)y=P(x)$ defined over a quadratic field
$K = \mathbb Q(\sqrt{D})$ that describe genus two curves with everywhere good reduction over $K$.

\section{An infinite family of genus two curves with everywhere good reduction}
\label{sec:infinite}

In this section we consider integral Weierstrass equations $E$ with $J_2(E)=0$. In \cite{CardonaQuer}, Mestre's algorithm was generalized to include genus 2 curves with automorphism group different from $\Z/2\Z$. These include genus two curves described by Weierstrass equations with the Igusa invariant $J_2=0$.  

In \cite[Theorem 3.11]{Shaska}, it was shown that if either $J_2=J_4=0$, or $J_2=J_6=0$, then the genus 2 curves that correspond to the moduli points $(J_2,J_4,J_6,J_{10})$ in $\mathcal{M}_2(\overline{\Q})$ are described by $1$-parameter families of Weierstrass equations. In the following theorem we consider the $1$-parameter family corresponding to $J_2=J_4=0$, and we show that for infinitely many rational values of the parameter these Weierstrass equations describe genus 2 curves with everywhere good reduction over real (complex) quadratic extensions of $\Q$. 

\begin{Theorem}
\label{thm:infinite}
Let $E_{j(t)}$ be the Weierstrass equation given as follows
{\footnotesize\begin{align*}
y^2 - x^2y  & =
(40960000000000 - 14580000000 j^5 + 531441 j^{10})\; x^6  -  (86400000000 j^2  - 39366000 j^7 )\;x^5\\&
- (1/4 + 48600000 j^4  + 177147 j^9 /4)\;x^4 + (320000000 j + 291600 j^6 )\;x^3 - 1350000 j^3 \;x^2 \\&
+ (1600000  + 243 j^5 )\;x-500 j^2
\end{align*} }
where $j(t)=20t-3$, $t$ is an integer. Then one has that $E_{j(t)}$ describes a genus 2 curve with everywhere good reduction over either of the extensions $K=\Q\left(\sqrt{\pm(-3200000 + 729 j(t)^5)}\right)$, for every integer $t$.
\end{Theorem}

\begin{Proof}
We take $J_{10}=\pm1$ in the parametric Weierstrass equation with Igusa invariants $J_2=J_4=0$ described in \cite[Theorem 3.11]{Shaska} (iv). We then minimize the Weierstrass equation at the primes $3,5$ to obtain the following $\Q(j)$-equivalent Weierstrass equation
{\footnotesize\begin{align*}
E'_j: y^2 & = 40960000000000 - 14580000000 j^5 + 531441 j^{10} +(- 172800000000 j^2  +
 78732000 j^7)\, x \\&
- (194400000 j^4 + 177147 j^9)\, x^2  + (2560000000 j  + 2332800 j^6)\, x^3 - 21600000 j^3 x^4 \\&
+  (51200000  + 7776 j^5)\, x^5 - 32000 j^2 x^6.
\end{align*}}
One notices that $J_{10}(E'_j)=-2^{30}(-3200000 + 729 j^5)^{15}$. 
For any integer value of $j$, it is clear that $E'_{j}$ is reduced over $\Q_p$ for any prime $p>5$. 

Taking $j\equiv -3$ mod $20$ assures that the primes $3$ and $5$ do not divide $J_{10}(E'_j)$, and moreover $2\nmid (-3200000 + 729 j^5)$. Setting $j=20t-3$, we use the transformation $x\mapsto 1/2x$, $y\mapsto(8y -4 x^2)/(2 x)$ to obtain the equation $E_{j(t)}:\;y^2 +Q(x) y=P(x)$, where $Q(x)=-x^2$, and 
{\footnotesize\begin{align*} 
P(x) & =
(40960000000000 - 14580000000 j^5 + 531441 j^{10})\; x^6  -  (86400000000 j^2  - 39366000 j^7 )\;x^5\\ &
-  (1/4 + 48600000 j^4  + 177147 j^9 /4)\;x^4 +  (320000000 j + 291600 j^6 )\;x^3 - 1350000 j^3 \;x^2 \\&
+ (1600000  + 243 j^5 )\;x-500 j^2.
\end{align*} }
One sees that $J_{10}(E_{j(t)})$ is given by $D(t):=-(-3200000 + 729 (20t-3)^5)^{15}$, hence $\gcd(J_{10}(E_{j(t)}),60)=1$ for any integer value of $t$. 

 One recalls that if $p^k$, $k\ge2$, divides $D(t)$ for some integer $t$, then $E_{j(t)}$ is not minimal over $\Q_p$, since otherwise the corresponding curve does not have potential good reduction over $\Q_p$, \cite[Remarque 13]{Liu}. Since the leading coefficient $\ell(t)$ of the sextic polynomial $4P(x)+Q(x)^2$ can be expressed as $\ell(t)=(-16800000 + 729 j(t)^5) (-3200000 + 729 j(t)^5)-2^{18}\times 5^{11}$, it follows that a common prime divisor of $D(t)$ and $\ell(t)$ should be one of the primes $2$ or $5$, yet our choice of $j(t)$ implies that $D(t)$ is odd and $5\nmid D(t)$. Therefore, any prime divisor of $D(t)$ does not divide $\ell(t)$. In addition, according to Corollary \ref{cor:irreducible}, one knows that $\overline{\omega}$ is non-ramified. The latter two facts imply that if $E_{j(t)}$ does not describe a curve with good reduction over $\Q_p$, then the minimal field $F$ over which good reduction is attained is a quadratic extension of $\Q_p$, \cite[Th\'eor\`eme 1 (a)]{Liu3}. 
 
Using Remark \ref{rem:Fequivalentnonreduced} and Remark \ref{rem:reduced}, one may assume without loss of generality that $E_{j(t)}$ is not $F$-equivalent to a non-reduced Weierstrass equation $E'_{j(t)}$ with $\nu_F(J_{10}(E_{j(t)}))=\nu_F(J_{10}(E'_{j(t)}))$. Therefore, Theorem \ref{thm:local} yields that $\nu_p(J_{10}(E_{j(t)}))$ is either $0$ or $15$. It follows that for every integer value of $t$, $E_{j(t)}$ satisfies the assumptions of Theorem \ref{thm1} for all prime divisors of $J_{10}(E_{j(t)})$. In addition, since $J_{10}(E_{j(t)})= D(t)=-(-3200000 + 729 (20t-3)^5)^{15}$, one has that $|-3200000 + 729 (20t-3)^5|=p_1\cdots p_n $, where $p_i$'s are the primes for which $\nu_p(J_{10}(E_{j(t)}))=15$. In view of Theorem \ref{thm1}, one sees that $E_{j(t)}$ describes a curve with everywhere good reduction over $\Q\left(\sqrt{\pm(-3200000 + 729 j(t)^5)}\right)$.
\end{Proof}
In the following example, we work out the details of how to use Theorem \ref{thm:infinite} to produce genus 2 curves that attain everywhere good reduction over both imaginary and real quadratic fields. 
\begin{Example}
In Theorem \ref{thm:infinite}, we consider the case when $t=0$. We set $C$ to be the curve described by the Weierstrass equation $E_{-3}$ over $\mathbb Q$. After the base change to $K^+=\Q(\sqrt{109\cdot 30983})$, the curve $C$ attains everywhere good reduction. In this case, $C$ is described by the following integral Weierstrass equation over $K^+$
{\footnotesize\begin{align*}y^2+a\;x^2 y&=692147988184705983745437384588460 + 
 1430511020616757602335438237 \cdot a\; x + 
 4160276333344024970164868400 \;x^2\\ &+ 
 1910741390971479497720 \cdot a\; x^3 +1667071152916855272634\; x^4 + 229696894990560 \cdot a\; x^5+ 
 44534321059609 \; x^6\end{align*}}
 where $a=\sqrt{109\cdot 30983}$. One can use \Magma to see that the discriminant of the latter Weierstrass equation is $1$.  

 In addition, after a base change to $K^-=\Q(\sqrt{-109\cdot 30983})$, the curve $C$ attains everywhere good reduction. If $b=\sqrt{-109\cdot 30983}$, the curve is described by the following integral Weierstrass equation over $K^{-}$
 {\footnotesize\begin{align*}y^2+x^2 y&= -692147988184705983745437384588460 - 
 1430511020616757602335438237\cdot b\; x + 
 4160276333344024970164868400\; x^2\\& + 
 1910741390971479497720\cdot b\; x^3 -1667071152916856116921\;x^4 - 229696894990560\cdot b\; x^5+ 
 44534321059609\;x^6
 \end{align*}}
 with discriminant $-1$.
\end{Example}

\begin{Corollary}
\label{cor:simple} Let $J_{j(t)}$, $t\in\Z$, be the Jacobian of the genus two curve described by $E_{j(t)}$ in Theorem \ref{thm:infinite}. Then $\End(J_{j(t)})=\Z$ for infinitely many integer values $t$.  In particular, $J_{j(t)}$ is an absolutely simple abelian variety. 
\end{Corollary}
\begin{Proof}
One may check using \Magma that if we set $t=1$ in $E_{j(t)}:y^2+Q_t(x)y=P_t(x)$, the Galois group of the polynomial $4P_1(x)+Q_1(x)^2$ is the symmetric group on six elements; $S_6$. Therefore, the polynomial $4P_t(x)+Q_t(x)^2$ must have generic Galois group $S_6$, in particular, there are infinitely many values of $t$ for which $4P_t(x)+Q_t(x)^2$ must have Galois group $S_6$. According to \cite[Theorem 2.1]{Zarhin}, for any curve $C$ described by an equation of the form $ y^2=f(x)$,  $n=\deg f \ge 5$, over a field $K$ with characteristic $0$, if $f$ has Galois group
$S_n$ or $A_n$ then $\End(\Jac(C)) = \mathbb Z$.  
\end{Proof}

\section{Final remarks and comments}

  (i) Let us mention the following "Open Question" stated in  \cite[p. 1382]{Dembele}.

\bigskip

    {\it Find a real quadratic field $F$, and an abelian surface $A$
    defined over $F$ such that $A$ has trivial conductor with
    $\End_F(A) = \mathbb Z$.}

\bigskip

We gave an answer to this question in Corollary \ref{cor:simple} using Jacobian varieties of genus two curves described by Weierstrass equations $E$ for which $J_2(E)=J_4(E)=0$, hence the automorphism group of the curves contain a nontrivial automorphism different from the hyperelliptic involution. Below we give an answer to the same question using the Jacobian of a genus two curve whose automorphism group is isomorphic to $\Z_2$,  in particular it is described by a Weierstrass equation $E$ with $J_2(E)\ne 0$. 

Consider the curve $C$ described by the Weierstrass equation $E: y^2=f(x)$ over the field $K = \mathbb Q(\sqrt{31 \cdot 47 \cdot 107 \cdot 269})$,
where

{\footnotesize\begin{align*}
f(x)&= 1/4 \cdot  (66861242074781803311769723545637855624989525807865019 \ -
 212180300372212133803172192824208735245404632954 \cdot a \; x  \\& +
11765761679797475148212972485522051834951241828993 \; x^2 -
8297325990561138463673721463315177816731190 \cdot  a \; x^3  \\& +
138030288410295455310019117974278390625714666 \; x^4 -
29202077961801148532323213140299077108 \cdot a \; x^5  \\& +
  107953688207841230963812605490049186473 \; x^6)
 \end{align*}}
and $a = \sqrt{31 \cdot 47 \cdot 107 \cdot 269}$. The Weierstrass equation $E$ is globally minimal except at the primes lying above $2$. Moreover, $E$ has Igusa invariants
$(J_2(E),J_4(E),J_6(E),J_{10}(E))=(1,8,-4,1)$.
Using \Magma, we may check that $C$ has everywhere good reduction over $K$ and the polynomial
$f(x)$ has the symmetric group $S_6$ as its Galois group. Now, in view of \cite[Theorem 2.1]{Zarhin},  we see that $\End(\Jac(C)) = \mathbb Z$.

(ii) Demb\'el\'e and Kumar in \cite{DembeleKumar} restricted the fields over which everywhere good reduction is attained to be real quadratic fields with narrow class number $= 1$.  We have examples of curves with  everywhere good reduction over
$K = \mathbb Q(\sqrt{D})$  and whose narrow class number is strictly larger than $1$.  One such example is the curve $y^2=f(x)$ over
$K = \mathbb Q(\sqrt{3244322821})$ with  $(J_2,J_4,J_6,J_{10})=(1, 6, 20,  1 )$ where 

{\tiny\begin{align*} 
f(x)&= 1/4 \cdot (56603909596961676599859889534368406320800421600598321259089511276547396625161821435597444
\\& -
   659996593966528308096994020555497202477042400305477225967492795853797578598449848 \cdot a \; x \\& +
   10402788757353302209521712779904492235857842873743815936622007573635398892236426913 \; x^2  \\& -
   26954577662737621340558324771176695069369536668994228423683609591814811006 \cdot a \; x^3 \\&+
   127456465699588438544545940903163923096316718284638517203001087373472362051 \; x^4 \\& -
   99075409972806731283591723191139320859584248134216074085707760902 \cdot a \; x^5 \\& +
   104107664637185788902792490762217692317930964923924351034101286541 \;x^6)
\end{align*}}
and $a = \sqrt{3244322821}$.
One easily checks that $h(3244322821) = 3$, hence the narrow class number
is strictly larger than $ 1$ (it is a multiple of the usual class number; more precisely it differs by a power of $2$).

(iii) In \cite[Table 7]{DembeleKumar}, a genus 2 curve $C$ with everywhere good reduction over $\Q(\sqrt{461})$ is given. 
The following example describes a genus 2 curve with everywhere good reduction over $\Q(\sqrt{461})$. 
The Jacobian of this curve is isogenous to the Jacobian of $C$. For $(J_2,J_4,J_6,J_{10})=(-36, -18, -4,  1 )$,
our algorithm gives the following minimal Weierstrass equation $y^2+Q(x)y=P(x)$ over $\Q(\sqrt{461})$ where

{\tiny\begin{align*}
Q(x)&=1,\\
P(x)&=12130483888208271296798313727552769926465009660 -
 598631135278559067562986005319927574362349047\cdot a\; x\\& +
 5674533451428904589539319523651961900875388818\; x^2 -
 62229859794374632967785456230016040735670239 \cdot a\; x^3\\& +
 176966449088848750932408731229660847277811408 \;x^4 -
 582211599034682482959237037970428723757123 \cdot a\; x^5 \\&+
 367925978657389359077234026398703090110814\; x^6
\end{align*}}
and $a = \sqrt{461}$.

We moreover obtain a curve with everywhere good reduction over $\Q(\sqrt{-461})$ described by the Weierstrass equation $y^2+Q^-(x)y=P^-(x)$ where
 {\tiny\begin{align*}
Q^-(x)&=a,\\
P^-(x)&=-12130483888208271296798313727552769926465009545 +
 598631135278559067562986005319927574362349047 \cdot a\;x \\&+
 5674533451428904589539319523651961900875388818\; x^2 -
 62229859794374632967785456230016040735670239 \cdot a\; x^3 \\&-
 176966449088848750932408731229660847277811408\; x^4 +
 582211599034682482959237037970428723757123 \cdot a\; x^5\\& +
 367925978657389359077234026398703090110814 \;x^6
\end{align*}}
and $a = \sqrt{-461}$.

\section*{Appendix}
\label{sec:Appendix}
 We fix $J_2=J_{10}= \pm1$. We consider $J_4$,  $1\le J_4\le 6$, and run the algorithm over the following values of $J_6$: $1\le J_6\le 50$. We were able to obtain
the following examples of globally minimal Weierstrass equations $y^2+Q(x)y=P(x)$ describing genus 2 curves with everywhere good reduction over a real (respectively, complex) quadratic field. In what follows $K_D=\Q(\sqrt{D})$, and $a_D=\sqrt{D}$.

For $(J_2,J_4,J_6,J_{10})=(1,3,41,1)$, we find the following minimal Weierstrass equation describing a genus 2 curve with everywhere good reduction over $K_D$, $D=283\cdot 983\cdot 292541$ with

{\tiny\begin{align*}
Q(x)&=-x^3 - x^2 - 2\;x - 1\\
P(x)&= 14493149125598449395929549754038474084310380241491931823193074055994900474640051957369500281928261566976890 \\&-(1+
 1665438444662855865052620278000056560398794643738231308609985628460383079986338377382156048929621\cdot a_D) \;x\\& +
 6489479481704636494295796125228374586785347202783533368433878328443971354643739735858926097025003\; x^2\\& - \frac{1}{2}(3+
 331431032874618735976084441445190757259603113874300807585890963582001508777869309138483\cdot a_D) \;x^3 \\&+
 193716095692334168253009328934954019939268327342175262407259127392365275308074763764493 \;x^4 \\&-\frac{1}{2}(1+
 2968043547749256931085921692531919016034618196595381301553295464182002510035\cdot a_D)\; x^5\\& +
 385505245693700062768252022867858762437647346167922980067063605138542968539\; x^6.
\end{align*}}

For $(J_2,J_4,J_6,J_{10})=(-1,3,-41,-1)$, we find the following minimal Weierstrass equation describing a genus 2 curve with everywhere good reduction over $K_D$, $D=-283\cdot 983\cdot 292541$ with

{\tiny\begin{align*}
Q(x)&= -x^3 - a_D\; x^2 - 2\;x - a_D\\
P(x)&=-14493149125598449395929549754038474084310380241491931823193074055994900474640051957369500281928241221554828\\& +
 1665438444662855865052620278000056560398794643738231308609985628460383079986338377382156048929620 \cdot a_D\; x \\&+
 6489479481704636494295796125228374586785347202783533368433878328443971354643739735858966787869128\; x^2\\& -
 165715516437309367988042220722595378629801556937150403792945481791000754388934654569243 \cdot a_D\; x^3\\& -
 193716095692334168253009328934954019939268327342175262407259127392365275308054418342433\; x^4\\& +
 1484021773874628465542960846265959508017309098297690650776647732091001255017 \cdot a_D\; x^5 \\&+
 385505245693700062768252022867858762437647346167922980067063605138542968539\; x^6.
\end{align*}}

For $(J_2,J_4,J_6,J_{10})=(1,4,4,1)$, we find the following minimal Weierstrass equation describing a genus 2 curve with
 everywhere good reduction over $K_D$, $D=26404561$ with

{\tiny\begin{align*}
Q(x)&=-x^3 - x - 1\\
P(x)&=3980068297155832423457786128424974929592052330976661679652318339872233\\
& -\frac{1}{2}(1+
 2615260210718839277044005360804018934630474400152774950170473767\cdot a_D) \;x\\
& +
 4726579929822021780213733516738163820518783843378859759003500803\; x^2\\
&-\frac{1}{2}(1+
 690174444833251297883378495386018560180595136487572674839 \cdot a_D)\;x^3 \\& +
 374207277614825042495998202528157743572394326529883851355\; x^4 \\&-
 8196253017457146326678417525184726280926437083474\cdot a_D\; x^5 \\& + 1975086830915274164097017785890796271744368850030 \;x^6.
  \end{align*}}

For  $(J_2,J_4,J_6,J_{10})=(-1,4,-4,-1)$, we find the following minimal Weierstrass equation describing a genus 2 curve
with  everywhere good reduction over $K_D$, $D=-26404561$ with

{\tiny\begin{align*}
Q(x)&=-x^3 - x -a_D\\
P(x)&=-3980068297155832423457786128424974929592052330976661679652318333271093\\& + 1307630105359419638522002680402009467315237200076387475085236883 \cdot a_D
\; x \\&+
 4726579929822021780213733516738163820518783843378859759003500803\; x^2\\& - 345087222416625648941689247693009280090297568243786337420\cdot a_D \; x^3 \\&-
 374207277614825042495998202528157743572394326529883851356\; x^4  \\&+
 8196253017457146326678417525184726280926437083474 \cdot a_D\;  x^5\\& + 1975086830915274164097017785890796271744368850030\; x^6.
\end{align*}}

For $(J_2,J_4,J_6,J_{10})=(-1,6,40,-1)$, we find the following minimal Weierstrass equation describing a genus 2 curve
with  everywhere good reduction over $K_D$, $D=22639\cdot3787361$ with

{\tiny\begin{align*}
Q(x)&=(-x^3 + a_D\cdot x^2 - 2\;x - a_D) \\
P(x)&=-20102088440135614030470128986946346563555781945937611141050719826788916139991909714976694933 \\&-
 4617292887655574165961428503504867769633166770554856646846593036443400760381748059 \cdot a_D\; x \\&-
 37889281434253382202455676997526996644934137545830281733194190489890863460961881935\; x^2 \\&-
 1933971633574852387913620699071895736642779234966827401936437145373079117\cdot a_D\;x^3 \\&-
 4761023220312967350217124372495946438441280174280887209278673812202315630\; x^4\\& -
 72904659362732357860282616825152037719939877323542568131715001\;a_D\; x^5 \\&-
 39883473190506942475701006247163415226409169450474757057166290 \;x^6.
\end{align*}}

For $(J_2,J_4,J_6,J_{10})=(1,6,-40,1)$, we find the following minimal Weierstrass equation describing a genus 2 curve
with  everywhere good reduction over $K_D$, $D=-22639\cdot3787361$ with
{\tiny\begin{align*}
Q(x)&=-x^3 + x^2 - 2\;x - 1\\
P(x)&= 20102088440135614030470128986946346563555781945937611141050719826788916139991909693541178513 \\&+(-1 +
 4617292887655574165961428503504867769633166770554856646846593036443400760381748058\cdot a_D\;) x\\& -
 37889281434253382202455676997526996644934137545830281733194190489890863503832914774\; x^2\\& +\frac{1}{2}(1 -3867943267149704775827241398143791473285558469933654803872874290746158235\cdot a_D)\; x^3 \\&+
 4761023220312967350217124372495946438441280174280887209278673790766799208\; x^4\\& +\frac{1}{2}(1 +
 145809318725464715720565233650304075439879754647085136263430003\cdot a_D
  )\; x^5\\& -
 39883473190506942475701006247163415226409169450474757057166290\; x^6.
\end{align*}}

\end{document}